\documentstyle{article}
\begin{document}
\newcommand{\Endo}[1]{\mathop{\mbox{End}_{#1}}}
\newcommand{\R}{\mbox{\boldmath $R$}}
\newcommand{\C}{\mbox{\boldmath $C$}}
\newcommand{\bfZ}{\mbox{\boldmath Z}}
\newcommand{\qH}{\mathop {\mbox{\boldmath $H$}}\nolimits}
\newcommand{\bfQ}{\mbox{\boldmath Q}}
\newcommand{\bfN}{\mbox{\boldmath N}}
\newcommand{\bfP}{\mbox{\boldmath P}}
\newcommand{\bfT}{\mbox{\boldmath T}}
\newcommand{\sqdz}{\sqrt{}}
\newcommand{\sqdzb}{\sqrt{d\bar{z}}}
\newcommand{\hD}{{\cal{D}}^\frac{1}{2}}
\newcommand{\calD}{{\cal{D}}}
\newcommand{\imH}{\mathop{\rm{im}}(\qH)}
\newcommand{\idx}{\mathop {\rm ind}\nolimits}
\newcommand{\setB}[1]{\left\{ #1 \right\}}
\newcommand{\divg}{\mathop {\rm div}\nolimits}
\newcommand{\cS}{{\cal S}} 
\newcommand{\fref}[1]{Figure~\ref{#1}}
\newcommand{\cM}{{\cal M}}
\newcommand{\MF}{{\cal MF}}
\newcommand{\mref}[1]{{\rm (\ref{#1})}}
\newcommand{\cT}{{\cal T}}
\newcommand{\s}{\mathop {\rm R}\nolimits}
\newcommand{\Sg}{(\Sigma, g)}
\newcommand{\RS}{\R_{+}^{\cal S}}
\newcommand{\Dif}{\mathop{\rm Diff_0}}
\newcommand{\pone}{\psi_{1}}
\newcommand{\ptwo}{\psi_{2}}
\newcommand{\None}{N_{1}}
\newcommand{\Ntwo}{N_{2}}
\newcommand{\Sione}{\Sigma_{1}}
\newcommand{\Sitwo}{\Sigma_{2}}
\newcommand{\II}{\mathop {\rm I\!I}\nolimits}
\newcommand{\IIone}{\mathop {\rm I\!I}_{1}\nolimits}
\newcommand{\IItwo}{\mathop {\rm I\!I}_{2}\nolimits}
\newcommand{\Sig}{(\Sigma, g)}

\newcommand{\n}{\noindent}
\newtheorem{theorem}{Theorem}
\newtheorem{assertion}{Assertion}
\newtheorem{proposition}{Proposition}
\newtheorem{remark}{Remark}
\newtheorem{lemma}{Lemma}
\newtheorem{definition}{Definition}
\newtheorem{claim}{Claim}[section]
\newtheorem{corollary}{Corollary}[section]
\newtheorem{coro}{Corollary}[section]
\newtheorem{observation}{Observation}[section]
\newtheorem{obs}{Observation}[section]
\newtheorem{conjecture}{Conjecture}[section]
\newtheorem{question}{Question}[section]

\title{Prescribing Mean Curvature: Existence and Uniqueness Problems}
\author{G. Kamberov \\ Department of Mathematics \\
Washington University, St. Louis, MO \\
e-mail: kamberov@math.wustl.edu}
\date{\  }

\maketitle

\section*{Introduction} 
The results in this paper are part of a program to understand 
to what extent one can use mean curvature data to 
determine the shape of a surface in $\R^3$. 

The first section of the paper presents results on a
problem suggested by Bonnet. A generic immersion is
uniquely determined up to a rigid motion by its first fundamental form
and its mean curvature function, but there are some exceptions, for
example most constant mean curvature immersions.  Bonnet's problem is
to classify and study all exceptional immersions.  Here I concentrate on
the study of Bonnet's problem for immersions with umbilic points and
immersions of closed surfaces.  

The second section contains an outline of an existence theory for conformal
immersions (\cite{GANG, franz, ulrich, saopaulo}) along with its 
immediate applications to Bonnet's problem.  The central idea in
the theory is to define square roots of basic geometric
objects, for example, area elements, differentials of maps, and then
determine the equations satisfied by the square root of the
differential of a conformal immersion.  From an analytic point of view
the main advantage of the approach is that the existence problem
for prescribing mean curvature data is reduced to a
first order elliptic system. The theory suggests a new paradigm: the
data used to determine the immersion are, the conformal structure, the
regular homotopy class of the immersion, and the mean curvature
half-density, that is, the mean curvature appropriately weighted by the
metric.  These data determine the immersion via a generalized
Weierstrass-Kenmotsu formula.

Robert Bryant introduced me to Bonnet's problem. I thank him and Ulrich
Pinkall, Franz Pedit, Peter Norman, and Dennis DeTurck for many discussions, advice, and support. 

\section{Bonnet's Uniqueness Problem}
Two isometric immersions $F_1$ and $F_2$, of a given Riemannian
surface $(M, g)$ into $\R^3$, are called congruent (denoted by 
$F_1 \sim F_2$) if they differ by a rigid motion. I am 
interested in the following questions: What conditions guarantee that
$H_1 = H_2$ implies $F_1 \sim F_2$? What can be said if $H_1 = H_2$
but $F_1$ and $F_2$ are not congruent? The isometric immersions $F_1$ and 
$F_2$ are called Bonnet mates if they are not congruent 
but their mean curvature functions, $H_1$ and $H_2$ are equal. 

Bonnet, Cartan, and Chern studied the existence and classification 
of Bonnet mates. Their works yield a complete local classification 
of the umbilic-free immersions which admit Bonnet mates (\cite{Bonnet},
\cite{Cartan-pr-curv}, \cite{Chern-Deform-Pr-Curv}). 
Some global results were obtained by Lawson
and Tribuzy, and by Ros (see \cite{trib_and_lawson} and 
\cite{ARoss}). 

The methods used in \cite{Bonnet}, \cite{Cartan-pr-curv}, and
\cite{Chern-Deform-Pr-Curv} are very powerful and beautiful but can not
be used to study non-constant mean curvature surfaces with umbilic
points. Until recently it was not known whether such 
Bonnet mates exist. In particular, it was not known
whether there are any Bonnet mates which are not included in the
Bonnet-Cartan-Chern classification. A general construction which
yields all immersions, with or without umbilic points, that admit
Bonnet mates was found in \cite{pipeka}. This construction made it
possible to prove that there exist infinitely many new Bonnet 
immersions, that is, Bonnet mates which have umbilic points and 
whose mean curvature is not constant in a neighborhood of the
umbilic points (\cite{saopaulo} and \cite{hqd}). 

This section contains results on the properties of the umbilic
points of Bonnet mates and on the the following rigidity
conjecture: {\it Let $(M,g)$ be a closed oriented Riemannian surface
  and let $H:M\rightarrow\R$ be an arbitrary non-constant function.
  Then up to rigid motions, there exists at most one isometric
  immersion $F:(M,g)\rightarrow \R^3$ with mean curvature function
  $H$}. 

Throughout this note $(M, g)$ is an oriented, connected Riemannian
surface; $F_1$, $F_2$ are isometric immersions of $(M, g)$ in
$\R^3$ ; $\II_i$, $N_i$, and $H_i$ denote the second
fundamental form, the field of unit normals, and
the mean curvature function of the immersion $F_i$.
The mean curvatures of two isometric immersions $F_1$ and $F_2$
coincide if and only if at every point $F_1$ and $F_2$ have the same 
principal curvatures, and so if $H_1 = H_2$ then $p$ is an
umbilic for one of the immersions precisely if it is umbilic for the 
other too.

It is convenient to introduce the shape distortion tensor, 
$D = \II_1 - \II_{2}$, associated with 
$F_1$ and $F_2$. Thus Bonnet's fundamental theorem of 
surface theory implies that $F_1 \sim F_2$ if and only if $D \equiv 0$. 
The Mainardi-Codazzi equations imply the following observation:
\begin{obs}\label{main-th-1} If the mean curvature functions 
of $F_1$ and $F_2$ coincide then the associated shape distortion
operator is trace-free and divergence free, i.e, 
its $(2,0)$ part, $D^{2,0}$, is holomorphic. 
\end{obs}
Observation \ref{main-th-1}  implies the following unique continuation 
property: 
\begin{coro} \label{unique-cont} Let $F_1: M\rightarrow \R^3$ and
  $F_2: M\rightarrow \R^3$ be isometric embeddings which have the same
  principal curvatures. If $F_1 \sim F_2$ on a nonempty open set then
  $F_1 \sim F_2$ on $M$.
\end{coro}
Analyzing the  shape distortion tensor yields the following theorem describing the properties of the umbilic points of a Bonnet immersion.  
\begin{theorem} \label{umbilics} Suppose that $F_1$ and $F_2$ are Bonnet 
mates then: 

{\bf A.}  $D_p = 0$ if and only if $p$ is umbilic; moreover, 
for every umbilic point $p$  we have 
\[ \idx_{F_1}(p) = \idx_{D^{2,0}}(p) =  \idx_{F_2}(p) < 0,\]
where $\idx_{F_1}(p)$ and $\idx_{F_2}(p)$ denote the index of 
$p$ with respect to $F_1$ and $F_2$ respectively and $\idx_{D^{2,0}}(p)$ 
is the index of $p$ with respect to the horizontal foliation of the 
quadratic differential $D^{2,0}$. (See  \cite{hopf} for the definitions of the 
different notions of index.) 

{\bf B.} Every umbilic point is a critical point of the mean and the 
Gauss curvatures. 
Furthermore, the trace-free part of $\II_{1}$ vanishes to a finite order 
at every umbilic point. 
\end{theorem}
According to Theorem \ref{umbilics}, if an isometric embedding, $F_1$,
admits a Bonnet mate then the net formed by the curvature lines of $F_1$ 
has the same local character, that is, the same type of 
singularities as the net
formed by the horizontal and vertical foliations of a holomorphic
quadratic differential. Next we consider what can be said about an
embedding whose foliations of curvature lines are precisely the
horizontal (vertical) foliations of a holomorphic quadratic
differential.
\begin{theorem} \label{folio} Let $M$ be a closed oriented 
  Riemannian surface  and let 
$F_{1}:(M,g)\rightarrow~\R^3$ be an isometric embedding, 
whose    
net of curvature lines is the net formed 
by the principal stretch foliations of a holomorphic quadratic 
differential. Then an isometric embedding
  $F_2 : (M,g)\rightarrow\R^3$ is congruent to $F_1$ if and
  only if $F_1\circ F_2^{-1}$ is orientation preserving and $H_1 =
  H_2$.
\end{theorem}
The proof of Theorem \ref{folio} also yields a theorem about immersions:
\begin{theorem} \label{foliimm} Let $M$ be a closed oriented 
  Riemannian surface  and let 
$F_{1}:(M,g)\rightarrow~\R^3$ be an isometric immersion, whose    
net of curvature lines is the net formed 
by the principal stretch foliations of a holomorphic quadratic 
differential.  If $F_1$ admits a Bonnet mate then 
its mean curvature is constant. 
\end{theorem}
In particular the shape of most globally isothermic embeddings of a
closed surface are determined by the mean curvature. For example, if $M$
has genus one and $F_{1}$ is globally isothermic, or if $F_1(M)$ is a
surface of revolution, then every other isometric embedding, $F_2$,
s.t., $H_1 = H_2$ must be congruent to $F_1$.

The unique continuation property implies the following theorem: 
\begin{theorem} \label{conformal} Let $F_1, F_2 : (M, g) \rightarrow
\R^3$ be two isometric embeddings then $F_1 \sim F_2$ if and only if the
following set of conditions are satisfied: (i) $F_1 \circ F_2^{-1}$ is
orientation preserving; (ii) $H_1 = H_2$; (iii) there exists a 
conformal map $\varphi : S^2 \rightarrow S^2$ such that $N_1 = \varphi\circ N_2$. 
\end{theorem}
The unique continuation also yields a 
generalization of the following classical result: {\it If $H_1 = H_2$ and 
$N_1 = N_2$ then $F_1 = \mbox{translation} \circ F_2$}, see 
\cite{Gardner1969}. 
\begin{theorem} Let $F_1, F_2 : (M, g) \rightarrow
\R^3$ be two isometric embeddings then $F_1 \sim F_2$ if and only if the
following set of conditions is satisfied: (i) $F_1 \circ F_2^{-1}$ is
orientation preserving; (ii) $H_1 = H_2$; (iii) $(N_1 - N_2)(M)$
lies in a half space.  
\end{theorem}
\section{Dirac Spinors and Conformal Immersions}
To obtain deeper rigidity results and to try to find and classify the
possible counterexamples to the rigidity conjecture one needs a new
method for constructing surfaces.
Such a theory is outlined here and then used to investigate 
the following modified rigidity conjecture: \\
{\it Let $(M,g)$ be a closed and oriented Riemannian surface. Given a function $H$ on $(M, g)$, and a regular homotopy class,
  ${\cal F}$, of immersions of $M$ into $\R^3$, then up to rigid
  motions there exists at most one isometric
  immersion $F\in{\cal F}$ with mean curvature function $H$}. \\
Note that if the modified conjecture is false then the rigidity
conjecture from the previous section is also false.

Instead of a new method for constructing surfaces one could try to use
the Weierstrass-Kenmotsu representation \cite{Kenmotsu79}.
To apply this representation one must solve Kenmotsu's system of differential equations satisfied 
by the differential $dF$ of a conformal immersion  $F$ with prescribed mean
curvature function $H$. The main difficulty is
that the system is non-homogeneous, second order, and nonlinear if $H$
is non-constant, and, in addition, the system is degenerate at points at
which $H=0$.  

The Weierstrass representation of minimal surfaces was reformulated in terms
of spinors in \cite{Sullivan89}.  These ideas were used in
\cite{KusnerSchmidt}, where the authors also indicated that there should be a
theory for general, not necessarily minimal, surfaces.  Indeed, such theories
were developed in \cite{Bobenko}, and later in \cite{kon} (see also
\cite{Taimanov}), and \cite{Jorg}. During the academic year 1995-1996 the
GANG seminar at the University of Massachusetts set out to investigate the
role of spinors in the geometry of immersed surfaces and to develop a general
theory of the spinor representation of surfaces.  The goal was to design an
efficient and useful calculus, and to give a transparent explanation of the
objects involved in the representation \cite{ulrich, saopaulo, GANG}.

From now on $M$ always denotes a Riemann surface, that is, an oriented
surface with a chosen conformal structure. In particular, $M$ comes equipped
with a maximal holomorphic atlas $\displaystyle{\{(U_\alpha, z_\alpha)\}}$.
The conformal structure allows us to make sense of non-negative two-forms,
and their square roots, the non-negative and the non-positive half-densities.
Indeed, the square roots of a nonnegative two-form $wdx_\alpha \wedge
dy_\alpha$ defined by $\pm\sqrt{w dx_\alpha \wedge dy_\alpha } = \pm\sqrt{w}
|dz_\alpha|$ are sections in the fiber bundle $\hD$ of half-densities.  Thus
the square root of a non-negative two-form, in particular an area element is
a half-density and vice versa the square of a half density is
a two-form. Half-densities are necessary to introduce a conformally invariant
notion of surface tension:
\begin{definition} Let $F$ be a conformal immersion of $M$ into $R^3$ 
inducing the area element  $dA$ by pulling back the 
Euclidean metric by $F$, and let $H$ be the mean curvature function of the 
immersion $F$. The half 
density $H\sqrt{dA}$ is called the mean curvature half density of $F$.
\end{definition}
A conformal immersion, $F$, of a Riemann surface into 
$\R^3$ defines a spin structure on the Riemann surface. The spin
structure characterizes uniquely the regular homotopy class of the
immersion. See \cite{Ulrich85} and \cite{Sullivan89}.  
The square root, $\sqrt{dF}$,
of the differential of $F$ is a section in the associated spinor
bundle, $\Sigma$.

For the rest of this paper we identify $\R^{3}$ with the imaginary 
quaternions $\imH$. A spinor bundle on a Riemann surface is a quaternionic 
line bundle $\Sigma$ with transition functions 
$k_{\beta\alpha}:U_\alpha\cap U_\beta \rightarrow \C$ such that 
$k_{\beta\alpha}^2 = \partial z_\alpha/\partial z_\beta$.  
A choice of a spinor bundle is equivalent to choosing 
a square root bundle of the bundle of conformal $\R^{3}$-valued one 
forms on $M$. Indeed, the following theorem is known:
\begin{theorem}  \label{suff} A spinor bundle on a Riemann surface 
$M$ is a quaternionic line bundle $\Sigma$ on $M$ with a chosen 
endomorphism $J\in \Endo{\qH}(\Sigma)$ and a nontrivial 
quaternionic-hermitian,  fiber-preserving pairing 
\[(\cdot, \cdot): \Sigma \times \Sigma \rightarrow T^{*}M\otimes\qH,  \]
so that $J^{2} = -1$, 
and for every two spinors $\psi, \phi \in \Sigma$ based at the same
point $p$ we have 
\begin{equation} \label{jconfprop}
(\phi, \psi)(JX) = (J\phi, \psi)(X) =(\phi, J\psi)(X)
\end{equation}
for every $X\in T_{p}M$. Here $JX$ denotes the action of the complex 
structure on the vector $X$. 
\end{theorem}
Note that for every spinor $\psi\in\Sigma$ the form  $(\psi,\psi)=\omega$ is 
imaginary quaternionic valued. Furthermore, $(\psi,\psi)$ is a conformal $\R^3 = \imH$ 
valued one-form on $M$. 
The spinor $\psi$ is interpreted as the square
root of $\omega$. Every spinor $\psi$ defines a non-negative 
half-density $|\psi|^{2} := |\!|(\psi,\psi)|\!|$,  
where $|\!|\cdot|\!|$ is the Euclidean norm in $\R^3$. The half-density  
$|\psi|^{2}$ measures the relative dilation associated with the 
conformal form $(\psi,\psi)$. 
A choice of a spinor bundle is equivalent to a
choice of a spin structure, that is a holomorphic square root of the
canonical bundle $T^{(1,0)}M^{*}$ of $M$.  Indeed, given a 
spinor bundle $\Sigma$, define the complex line bundle of 
{\em positive spinors}  by $\Sigma_{+} \stackrel{{\rm def}}{=} 
\setB{\sigma\in \Sigma \, | \, J\sigma = \sigma i}$. For every 
positive spinor $\sigma\in \Sigma_{+}$ we define the one-from 
$\sigma^{2} \stackrel{{\rm def}}{=}-k(\sigma, \sigma)$. From 
\mref{jconfprop} it follows that 
$\sigma^{2}(JX)=i\sigma^{2}(X) = \sigma^{2}(X)i$ for every vector 
$X$. Thus for every positive spinor $\sigma$, the one-form
$\sigma^{2}$ is complex valued and belongs to the canonical bundle 
$T^{(1,0)}M^{*}$ of $M$. The map sending $\sigma\in \Sigma_+$ to 
$\sigma^{2}$  
is a quadratic map from $\Sigma_{+}$ onto the canonical bundle of $M$.
Thus every $\psi\in \Sigma_+$ defines a conformal one form 
$\psi^{2}\in T^{(1,0)}M^*$ and a non-negative half density 
$|\psi|^2$. 

There is a canonical bi-additive fiber-wise pairing $\Sigma \times (T^*M
\otimes_{R}\qH) \rightarrow D^{\frac{1}{2}} \otimes_{R} \Sigma$ which 
defines  conformal Clifford product between spinors and 
$\qH$-valued one-forms. 
The {\it conformal Dirac operator} is the unique local linear operator $ D :
\Gamma(\Sigma) \rightarrow\Gamma(\calD^{\frac{1}{2}}\otimes\Sigma) $
satisfying the Leibniz rule and such that $D\sigma = 0$ if $\sigma$ is a
local section of $\Sigma_+$ whose square is a closed one-form, i.e., if
$d\sigma^2 = 0$.  (Compare with \cite{Hitchin1974} and \cite{Atyah71}.)
\begin{definition} A section,  $\psi\in \Gamma(\Sigma)$, is called 
  a Dirac spinor if $D\psi = U\otimes \psi$ for some $U\in \Gamma(\hD)$; we
  say that $\psi$ generates the half density $U$.
\end{definition}
Dirac spinors represent the square roots of the differentials of conformal
immersions:
\begin{theorem} \label{maintheorem}
  For every non-vanishing $\psi\in \Gamma(\Sigma)$ the one-form $(\psi,
  \psi)$ is closed if and only if $\psi$ is a Dirac spinor. If $M$ is simply
  connected and $\psi$ is a non-vanishing Dirac spinor, then
\begin{equation} \label{weirstrass} 
F = \int (\psi, \psi)
\end{equation}
is a conformal immersion with mean curvature half density equal to the
half-density generated by $\psi$.  Vice versa, given an oriented
surface $M$,
then every immersion $F$ of $M$ into $\R^3$ defines a unique up to
isomorphism spinor bundle
$\Sigma$ on $M$, the surface $M$ is equipped with the pull-back conformal structure, and precisely two Dirac spinors, $\psi$, and $-\psi$ such
that $(\pm\psi, \pm\psi) = dF$. The half-density generated by $\pm\psi$
equals the mean curvature half-density of $F$. 
\end{theorem}
Note that in fact $\Sigma = M\times\qH$ and after identifying $\R^3$ with 
$\imH$ we have $(\psi, \phi) = \overline{\psi}dF\phi$, where
$\psi$ and $\phi$ are spinors based at the same point on $M$.
I will denote by $[\Sigma]$ the regular homotopy class of an immersion
inducing the spinor bundle $\Sigma$. 

The standard elliptic theory and the representation \mref{weirstrass}
imply the following local  existence result: {\it Every half-density can be realized
  locally as the mean-curvature half-density of a conformal immersion}.
Global existence and the Dirichlet problem are
discussed in \cite{diracexists} and \cite{dirichlet}.

The outlined theory provides new insight into the Bonnet's 
problem. In particular, in the structure of the space of pairs
of Bonnet mates on a closed surface $M$. Given a metric $g$ on $M$ and a 
function $H\neq \mbox{const}$ then up to
rigid motions there exists at most one Bonnet pair $f_{\pm}$
of isometric immersions of $(M,g)$ into $\R^3$ with mean curvatures  
 $H_{\pm}=H$ (\cite{trib_and_lawson}). I am interested in the space 
of {\bf geometrically distinct} pairs of
Bonnet mates on a given Riemann surface $M$ and within the 
same regular homotopy class. The space of geometrically distinct pairs
of 
Bonnet mates is the space of Bonnet pairs modulo 
the natural gauge group ${\cal G} = {\cal E}(3) \times {\cal E}(3) \times \R^{+}$
acting on it. Here ${\cal E}(3)$ is the 
Euclidean group of rigid motions. Indeed, if we rotate one of the 
Bonnet mates in the pair $F_{\pm}$ we obtain another pair of Bonnet
mates which is geometrically identical with the original pair
$F_{\pm}$. Furthermore, for every positive 
number $r > 0$ the Bonnet pair $rF_{\pm}$ is simply the original pair 
$F_{\pm}$ observed at a different scale, thus the pairs $F_{\pm}$ and 
$rF_{\pm}$ are not geometrically distinct. The 
action of the group ${\cal G}$ 
preserves the conformal class, the mean curvature half density and the
regular homotopy class. Let ${\cal B}(M, \Sigma, U)$ be 
the moduli space of  geometrically distinct 
pairs of Bonnet mates defined on the Riemann surface $M$ which belong
to the regular homotopy class associated with the spinor bundle
$\Sigma$, and inducing a given potential $U$. 

\begin{theorem} \cite{diracexists} \label{ker2} For every half-density $U$, 
regular homotopy class $[\Sigma]$ on the closed
  Riemann surface $M$ the moduli space ${\cal B}(M, \Sigma, U)$ of
  geometrically distinct Bonnet mates is either empty or it is a
  disjoint union of isolated components. Each nonempty component is either a
  point, a line, or a four dimensional ball.
\end{theorem} 

Note that $\dim_{\qH}\ker(D - U)$ is finite if $M$ is closed. 
If $\dim_{\qH}\ker(D - U)=1$, then ${\cal B}(M, \Sigma, U)$ is empty. Furthermore:

\begin{theorem} \cite{diracexists} \label{dim2} Suppose that $U$ is a
  half density on the Riemann surface $M$, if the spinor bundle $\Sigma$ is
  such that $\dim_{\qH}\ker(D - U) = 2$, then precisely one of the
  following alternatives holds: \\

{\bf(a.) } ${\cal B}(M, \Sigma, U) = \emptyset$. \\

{\bf (b.) } ${\cal B}(M, \Sigma, U)$ is a point. \\

{\bf (c.) } ${\cal B}(M, \Sigma, U)$ is homeomorphic to a line segment. \\

{\bf (d.) } ${\cal B}(M, \Sigma, U)$ is homeomorphic to a
four dimensional ball. \\

Moreover,  every conformal immersion in the 
regular homotopy class $[\Sigma]$ and inducing the mean curvature half density 
$U$ admits at most one, up to rigid motions, Bonnet mate. 
\end{theorem}

Theorem~\ref{dim2} provides an extension of
the results from \cite{trib_and_lawson} including constant mean curvature 
immersions. 

\begin{corollary} \label{l-tr-ext} Let $(M,g)$ be an oriented, closed
  Riemannian surface, and let $H$ be a function defined on
$M$, possibly constant. Every regular homotopy class
$[\Sigma]$ such that $\dim_{\qH}\ker(D - U) = 2$ contains at most 
two geometrically distinct isometric immersions of $(M, g)$ into 
$\R^3$ with mean curvature function $H$. 
\end{corollary}

\noindent {\bf Remark}  One can draw a parallel between Bonnet's rigidity 
conjecture and Pauli's exclusion principle by thinking of the Dirac spinor
$\psi$ of an immersion $f$ as a wave function. Indeed, the Dirac spinor
satisfies the equation $D\psi = U\psi$ and defines the (probability)
half-density $|\psi|^2 = \sqrt{dA}$, where $dA$ is the area element induced
on the surface by the immersion $f$.  Two non-congruent immersions $f_{\pm}$
are Bonnet mates if and only if they are in the same conformal class, induce
the same area element and the same mean curvature half density $U_{+} =
U_{-}$. Thus if we consider only Bonnet mates within the same regular
homotopy class, then two non-congruent immersions $f_{\pm}$ are Bonnet mates
if and only if their respective Dirac spinors $\psi_{\pm}$ generate the same
potential $U_{+} = U_{-}$ and the same half-density $|\psi_{+}|^2 =
|\psi_{-}|^2$.  Furthermore, it is natural to think that the Dirac spinors
$\psi$ and $\psi\alpha$, where $\alpha$ is a unit quaternion, represent the
wave functions of the same particle, that is the same quantum state, with
respect to two different Euclidean reference frames in the ambient space
$\R^3$. Hence the modified rigidity conjecture is equivalent to saying that
the quantum state is uniquely determined by the potential (energy) and the
position density.

\bibliographystyle{plain}

\end{document}